\documentclass{amsart}%[11pt]{article}

\newcommand{\Z}{\mathbb Z}
\newcommand{\W}{\mathbb W}
\newcommand{\R}{\mathbb R}

\usepackage{amsfonts,amssymb, amsthm}

\newtheorem{theorem}{Theorem}[section]

\newtheorem{lemma}[theorem]{Lemma}
\newtheorem{proposition}[theorem]{Proposition}

\theoremstyle{definition}

\input{epsf}

\usepackage{graphics}

\date{}

\begin{document}

%%%%%%%%%%%%%%%%%
%%%%%%%%%%%%%%%%%

\title[Type III moves and Fox colorings]{
A lower bound for the number of \\
Reidemeister moves of type III}

\author{J. Scott Carter} 
\address{Department of Mathematics, 
University of South Alabama, Mobile, AL 36688, U.S.A.} 
\email{carter@jaguar1.usouthal.edu} 

\author{Mohamed Elhamdadi} 
\address{Department of Mathematics, 
University of South Florida, Tampa, FL 33620, U.S.A.} 
\email{emohamed@math.usf.edu}

\author{Masahico Saito} 
\address{Department of Mathematics, 
University of South Florida, Tampa, FL 33620, U.S.A.} 
\email{saito@math.usf.edu} 

\author{Shin Satoh} 
\address{Graduate School of Science and Technology, 
Chiba University, Yayoi-cho 1-33, 
Inage-ku, Chiba, 263-8522, Japan 
(Department of Mathematics, University of South Florida, 
April 2003--March 2005)}
\email{satoh@math.s.chiba-u.ac.jp}

\dedicatory{Dedicated to Professor Louis H. Kauffman for his $60$th birthday}

\renewcommand{\thefootnote}{\fnsymbol{footnote}}


\begin{abstract}
We study the number of Reidemeister type III moves using Fox
$n$-colorings of knot diagrams.
\end{abstract}

\maketitle

%%%%%%%%%%%%%%%%%
%%%%%%%%%%%%%%%%%
\section{Introduction}
%%%%%%%%%%%%%%%%%
%%%%%%%%%%%%%%%%%

\begin{figure}[htb]
\begin{center}
\includegraphics{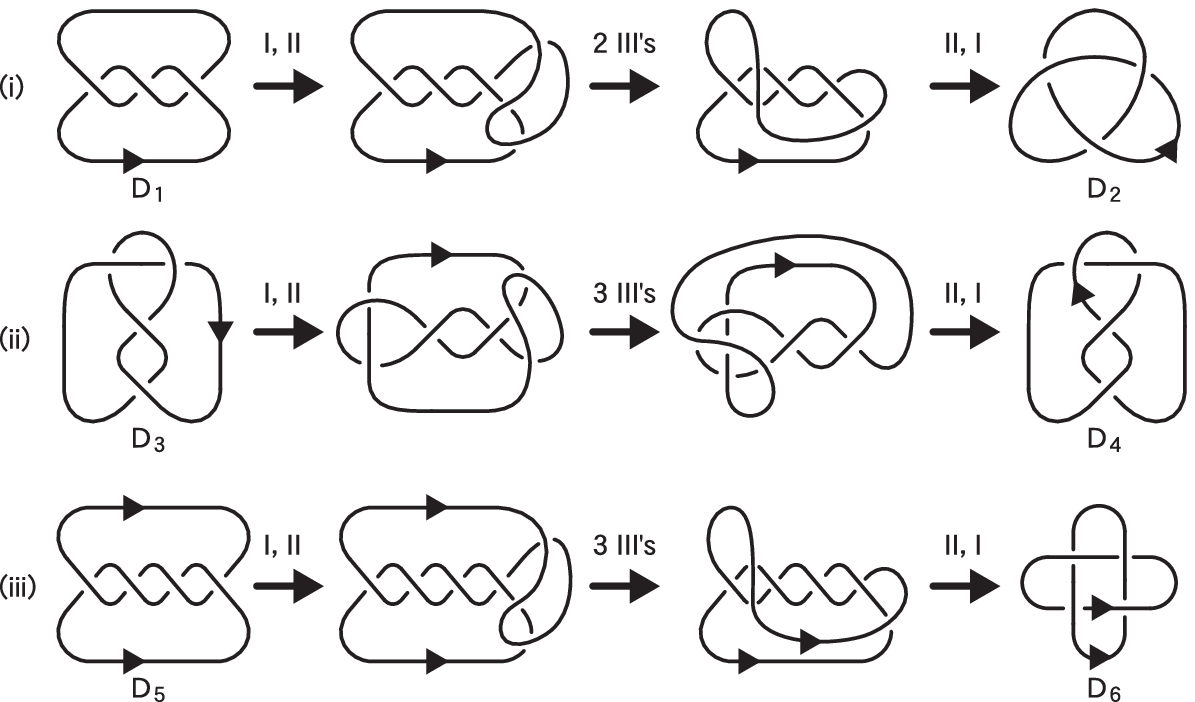}
\caption{}
\label{fig01}
\end{center}
\end{figure}

Any two oriented diagrams in the plane $\R^2$ 
representing the same knot or link 
are related by a finite sequence 
of Reidemeister moves of type I, II, and III
 and planar isotopies on the underlying graph.
For instance, 
Figure~\ref{fig01} illustrates 
three pairs of oriented diagrams representing 
the (i) trefoil knot, (ii) figure-eight knot, 
and (iii) $(2,4)$-torus link, respectively. 
The reader can construct 
such finite sequence of Reidemeister moves 
for each pair. 

%%%%%%%%%%%%%%%%%
%%%%%%%%%%%%%%%%%

%One of the topics of the paper is 
The main topic of the paper is %%%
the minimal number of Reidemeister moves of type III 
for all possible sequences 
relating two given oriented diagrams: 
For a pair of oriented diagrams $(D,D')$ of a knot or link, 
we denote by $\Omega_3(D,D')$ 
the minimal number of type III moves 
connecting $D$ and $D'$. 
Then we prove the following.

%%%%%%%%%%%%%%%%%
\begin{theorem}\label{thm11}
The oriented diagrams $D_i$ $(i=1,2,\dots,6)$ 
in {\rm Figure~\ref{fig01}} satisfy
\begin{itemize} 
\item[{\rm (i)}] $\Omega_3(D_1,D_2)=2$, 
\item[{\rm (ii)}] $\Omega_3(D_3,D_4)=3$, 
and 
\item[{\rm (iii)}] $\Omega_3(D_5,D_6)=3$. 
\end{itemize} 
\end{theorem} 
%%%%%%%%%%%%%%%%%

\begin{figure}[htb]
\begin{center}
\includegraphics{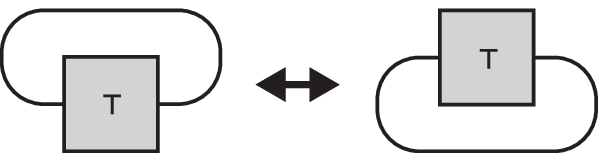}
\caption{}
\label{fig02}
\end{center}
\end{figure}

We remark that if a diagram is considered as lying 
on the $2$-sphere $S^2=\R^2\cup\{\infty\}$, 
then an outside arc of the diagram 
can be thrown through $\infty$, 
see Figure~\ref{fig02}. 
By allowing this move, 
the diagram $D_{2i-1}$ $(i=1,2,3)$ in Figure~\ref{fig01} 
can be deformed into $D_{2i}$ 
without any Reidemeister moves. 
%% add:
We also remark that the number of type III moves is related
to the minimal number of triple points in projections in $3$-space of 
embedded surfaces in $4$-space \cite{SatShi}. %% end add

%%%%%%%%%%%%%%%%%
%%%%%%%%%%%%%%%%%

This paper is organized as follows: 
In Section~\ref{sec2}, 
we review the Fox coloring for diagrams, 
which will be used to prove Theorem~\ref{thm11} 
in Section~\ref{sec3}.

%%%%%%%%%%%%%%%%%
%%%%%%%%%%%%%%%%%
\section{Preliminaries}\label{sec2} 
%%%%%%%%%%%%%%%%%
%%%%%%%%%%%%%%%%%

\begin{figure}[htb]
\begin{center}
\includegraphics{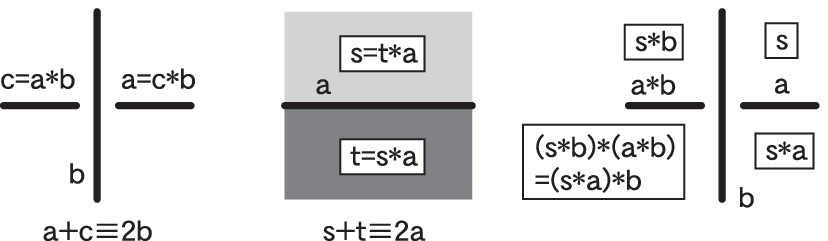}
\caption{}
\label{fig03}
\end{center}
\end{figure}

Let $D$ be a (possibly unoriented) link diagram in $\R^2$, 
which is an illustration of the projection of a link 
with small gaps at crossings 
to indicate over-under information. 
Thus $D$ is regarded as a disjoint union of arcs, 
and we denote by ${\rm Arc}(D)$ the set of the arcs of $D$. 
Let $n$ be a positive integer, 
and $\Z(n)=\{0,1,\dots,n-1\}$ 
the set of integers between $0$ and $n-1$ inclusive. 
Given a map $\varphi:{\rm Arc}(D)\rightarrow \Z(n)$, 
we call $\varphi(\alpha)$ the {\it color} of 
an arc $\alpha\in{\rm Arc}(D)$. 
We say that this map $\varphi$ is 
a {\it Fox $n$-coloring} 
\cite{Fox} %%% added
or simply 
{\it $n$-coloring} for $D$ 
if the equality $a+c\equiv 2b$ $({\rm mod}\ n)$
holds at each crossing of $D$, 
where $a$ and $c$ are the colors of 
the two under-arcs, 
and $b$ is the color of 
the over-arc at the crossing. 
See the right of Figure~\ref{fig03}, 
where we denote by $x*y$ the integer $k$ in $\Z(n)$ 
satisfying $k\equiv 2y-x$ (mod $n$) 
for $x,y\in\Z(n)$.

%%%%%%%%%%%%%%%%%
%%%%%%%%%%%%%%%%%

In this paper, 
we use  {\it extended} $n$-colorings
that were originally introduced in 
\cite{RS} (see also \cite{CKS1}). %% added
Let $D_*$ be the set of immersed circle(s) in $\R^2$ 
obtained from $D$ by ignoring crossing information, 
and ${\rm Region}(D)$ be 
the set of connected regions of the complement 
$\R^2\setminus D_*$. 
A map $\overline{\varphi}:
{\rm Arc}(D)\cup{\rm Region}(D)\rightarrow \Z(n)$ 
is called an extended $n$-coloring 
if (i) $\varphi=\overline{\varphi}|_{{\rm Arc}(D)}$ 
is the original Fox $n$-coloring, 
and (ii) the equality $s+t\equiv 2a$ $({\rm mod}\ n)$ 
holds at every point on $D$ except crossings, 
where $s$ and $t$ are the colors of 
the two regions around the point, 
and $a$ is the color of the arc 
on which the point lies. 
See the center of Figure~\ref{fig03}. 
We remark that there are no additional conditions 
around a crossing as shown in the right figure; 
indeed, it holds that 
$$\left\{
\begin{array}{l}
(s*b)*(a*b)\equiv 2(2b-a)-(2b-s)\equiv -2a+2b+s,\\
(s*a)*b\equiv 2b-(2a-s)\equiv -2a+2b+s 
\end{array}\right.$$
modulo $n$, 
which implies that $(s*b)*(a*b)=(s*a)*b$ in $\Z(n)$.

Given an $n$-coloring $\varphi$ for $D$ 
and an integer $s\in\Z(n)$, 
there is a 
unique extended $n$-coloring $\overline{\varphi}$ 
for $D$  such that the color of the outermost region 
in ${\rm Region}(D)$ is $s$. 
We denote it by $\overline{\varphi}=(\varphi,s)$. 
We say that an extended $n$-coloring 
$\overline{\varphi}=(\varphi,s)$ 
is {\it trivial} if $\varphi$ is a constant map, 
that is, all the arcs have the same color. 
Otherwise, $\overline{\varphi}$ is called {\it non-trivial}.

\begin{figure}[htb]
\begin{center}
\includegraphics{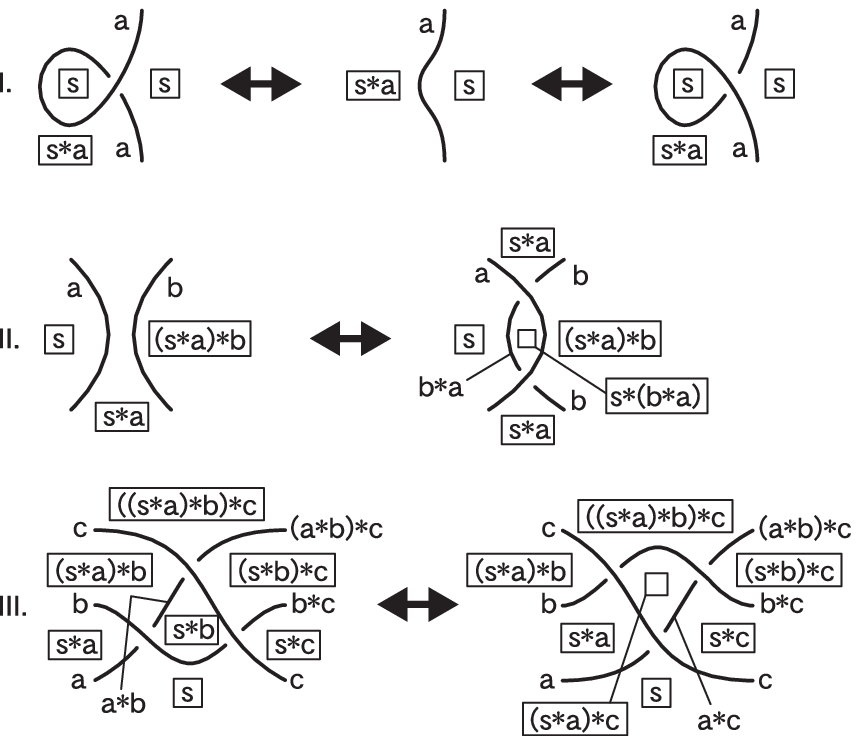}
\caption{}
\label{fig04}
\end{center}
\end{figure}

Assume that a diagram $D'$ is obtained from $D$ 
by a single Reidemeister move. 
If $D$ has an extended $n$-coloring $\overline{\varphi}$, 
then there is a unique  extended $n$-coloring 
$\overline{\varphi}'$ for $D'$  
such that these two colorings are coincident 
in the exterior of the disk where the Reidemeister move 
is performed. 
Figure~\ref{fig04} shows 
such Reidemeister moves between the pairs 
$(D,\overline{\varphi})$ and $(D',\overline{\varphi}')$. 
Note that if $\overline{\varphi}=(\varphi,s)$, 
then $\overline{\varphi}'=(\varphi',s)$ 
for some $n$-coloring $\varphi'$ for $D'$; 
that is, Reidemeister moves keep the color 
of the outermost region. 
Also, if $\overline{\varphi}$ is non-trivial, 
then so is $\overline{\varphi}'$. 

%%%%%%%%%%%%%%%%%
%%%%%%%%%%%%%%%%%

\begin{figure}[htb]
\begin{center}
\includegraphics{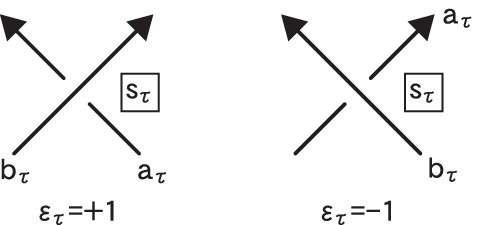}
\caption{}
\label{fig05}
\end{center}
\end{figure}

Let $f:\Z(n)^3\rightarrow \Z$ be 
a function from the set $\Z(n)^3=\Z(n)\times\Z(n)\times\Z(n)$ 
to $\Z$ which satisfies the condition 
$(\#)$ if $y=z$ then $f(x,y,z)=0$. 
For a pair $(D,\overline{\varphi})$ 
of an oriented diagram $D$ and 
its extended $n$-coloring $\overline{\varphi}$, 
we define an integer 
$\W_f(D, \overline{\varphi})$ 
as follows: 
Around a crossing $\tau$ of $D$, 
let $a_{\tau}\in\Z(n)$ be the color of the under-arc 
on the right side of the over-arc with respect 
to the orientation of the over-arc, 
let $b_{\tau}\in\Z(n)$ be the color of the over-arc, 
and let $s_{\tau}\in\Z(n)$ be the color of the region 
on the right side of the over- and under-arcs both 
with respect to their orientations. 
Also, let $\varepsilon_{\tau}\in\{+1,-1\}$ be 
the sign of the crossing $\tau$. 
The left and right of Figure~\ref{fig05} 
show such a triple $(s_{\tau}, a_{\tau}, b_{\tau})$ 
with $\varepsilon_{\tau}=+1$ and $-1$, 
respectively. 
Then we define 
$\W_f(D,\overline{\varphi})=
\sum_{\tau} \varepsilon_{\tau}
f(s_{\tau},a_{\tau},b_{\tau})\in\Z$, 
where the sum is taken over all the crossings $\tau$ of $D$.

%%%%%%%%%%%%%%%%%
%%%%%%%%%%%%%%%%%

Associated with a function $f:\Z(n)^3\rightarrow\Z$ 
which satisfies the property $(\#)$, 
we define the function 
$\delta f:\Z(n)^4\rightarrow\Z$ by 
\begin{eqnarray*}
(\delta f)(x,y,z,w) &=& 
f(x,z,w)-f(x,y,w)+f(x,y,z)\\
& & -f(x*y,z,w)+f(x*z,y*z,w)-f(x*w,y*w,z*w). 
\end{eqnarray*}
This definition is motivated from the coboundary operator 
of rack and quandle homology theories \cite{CJKLS,FRS0}. %% added

%%%%%%%%%%%%%%%%%
\begin{lemma}\label{lem21}
Assume that a pair $(D,\overline{\varphi})$ 
of an oriented diagram and its extended $n$-coloring 
is related to $(D',\overline{\varphi}')$ 
by a single Reidemeister move. 
\begin{itemize} 
\item[{\rm (i)}] 
If the move is of type {\rm I} or {\rm II}, 
then 
$\W_f(D,\overline{\varphi})=\W_f(D',\overline{\varphi}')$. 
\item[{\rm (ii)}] 
If the move is of type {\rm III}, 
then 
$\W_f(D,\overline{\varphi})-
\W_f(D', \overline{\varphi}')\in\pm{\rm Im}(\delta f)$. 
\end{itemize} 
\end{lemma} 
%%%%%%%%%%%%%%%%%

\begin{proof} 
(i) 
If the move is of type I 
in the first row of Figure~\ref{fig04}, 
we obtain 
$\W_f(D,\overline{\varphi})-
W_f(D', \overline{\varphi}')=\pm f(s,a,a)=0$ 
by the condition $(\#)$. 
Assume that the move is of type II 
in the second row of the figure. 
If the two arcs are both oriented downward, 
then it holds that 
$\W_f(D,\overline{\varphi})-
W_f(D', \overline{\varphi}')=
\pm\bigl[f(s,b*a,a)-f(s,b*a,a)\bigr]=0$. 
Other cases are similarly proved. 

(ii) 
We may assume that 
the three arcs in the bottom row 
are oriented from left to right. 
Then all the crossing are positive, 
and it holds that 
\begin{eqnarray*}
\W_f(D,\overline{\varphi})-
W_f(D', \overline{\varphi}') 
&=& \pm\bigl[
f(s,a,b)+f(s*b,a*b,c)+f(s,b,c)\\
& & -f(s*a,b,c)-f(s,a,c)-f(s*c,a*c,b*c) \bigr]\\
&=&
\pm(\delta f)(s,a,b,c),
\end{eqnarray*}
which belongs to the set $\pm{\rm Im}(\delta f)$. 
\end{proof}

Let $D$ be an oriented diagram, 
$f:\Z(n)^3\rightarrow \Z$ 
a function with the property $(\#)$, 
$s$ an integer in $\Z(n)$, 
and $m$ a non-negative integer. 
We define two finite sets of integers 
as follows: 

\begin{itemize} 
\item 
$\Phi_f(D,s)=
\bigl\{\W_f(D,\overline{\varphi})|\ 
\overline{\varphi}=(\varphi,s)
\mbox{ is a non-trivial $n$-coloring for } D\bigr\}$. 
\item
$\Delta_m(f)=\bigl\{
k_1+k_2+\dots+k_m|\ 
k_1,k_2,\dots,k_m\in\pm{\rm Im}(\delta f)\bigr\}$.
\end{itemize} 
Here we put $\Delta_0(f)=\{0\}$ 
for convenience.

%%%%%%%%%%%%%%%%%
\begin{proposition}\label{prop22}
Let $D$ and $D'$ be oriented diagrams 
of the same knot or link. 
If there is a non-trivial extended $n$-coloring 
$\overline{\varphi}=(\varphi,s)$ for $D$ 
and a function $f:\Z(n)^3\rightarrow\Z$ 
with the property $(\#)$ 
such that 
$$\bigl[\W_f(D,\overline{\varphi})-\Phi_f(D',s)\bigr]
\cap \Delta_i(f)=\emptyset$$
for any $i=0,1,\dots,m-1$, 
then 
$\Omega_3(D,D')\geq m$ holds, 
where $\W_f(D,\overline{\varphi})-\Phi_f(D',s)$ 
stands for the finite set of integers 
$\bigl\{\W_f(D,\overline{\varphi})-w|\ 
w\in\Phi_f(D',s)\bigr\}$. 
\end{proposition} 
%%%%%%%%%%%%%%%%%

\begin{proof} 
Put $\omega=\Omega_3(D,D')$. 
There is a finite sequence of Reidemeister moves 
between $D$ and $D'$, 
in which there are $\omega$ moves of type III. 
For the extended $n$-coloring 
$\overline{\varphi}'=(\varphi',s)$ for $D'$ 
associated with $\overline{\varphi}=(\varphi,s)$ for $D$, 
it follows that 
$\W(D,\overline{\varphi})-
\W(D',\overline{\varphi}')\in\Delta_{\omega}(f)$ 
by Lemma~\ref{lem21}. 
Hence, we have 
$\omega\geq m$. 
\end{proof}

%%%%%%%%%%%%%%%%%
%%%%%%%%%%%%%%%%%
\section{Proof of Theorem~$\ref{thm11}$}\label{sec3} 
%%%%%%%%%%%%%%%%%
%%%%%%%%%%%%%%%%%

\begin{figure}[htb]
\begin{center}
\includegraphics{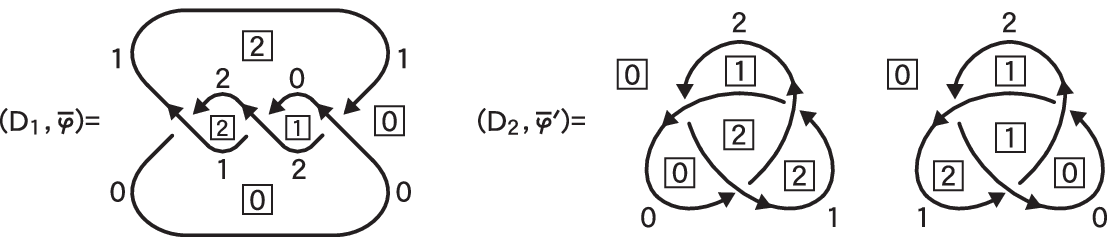}
\caption{}
\label{fig06}
\end{center}
\end{figure}

(i) We take a function 
$\Z(3)^3\rightarrow\Z$ 
defined by $f(x,y,z)=(x-y)(y-z)z$, 
which satisfies the condition $(\#)$. 
This definition of $f$ by a polynomial is motivated from \cite{Mochi}. %%% added
Let $\overline{\varphi}=(\varphi,0)$ 
be the extended $3$-coloring for $D_1$ 
as shown in the left of Figure~\ref{fig06}. 
For the left crossing, say $\tau$, 
we have the triple 
$(s_{\tau}, a_{\tau}, b_{\tau})=(2,2,1)$ by definition. 
Similarly, the center and right crossings 
have the triples 
$(2,0,2)$ and $(2,1,0)$, 
respectively. 
Since the signs of these crossings are all positive, 
the integer $\W_f(D_1;\overline{\varphi})$ 
is calculated by 
\begin{eqnarray*}
\W_f(D_1;\overline{\varphi}) &=& 
+f(2,2,1)+f(2,0,2)+f(2,1,0)\\
 &=& (2-2)(2-1)1+(2-0)(0-2)2+(2-1)(1-0)0
 =-8. 
 \end{eqnarray*}

%%%%%%%%%%%%%%%%%
%%%%%%%%%%%%%%%%%

%Any non-trivial extended $3$-coloring 
%$\overline{\varphi}'=(\varphi',0)$ 
%for the diagram $D_2$ has 
%two cases 
%(up to rotations of $D_2$) as 
%shown in the right of Figure~\ref{fig06}. 
There are two cases for non-trivial extended $3$-colorings
$\overline{\varphi}'=(\varphi',0)$ 
for the diagram $D_2$, that are shown in the right of Figure~\ref{fig06}. %%%
The integer $\W_f(D_2,\overline{\varphi}')$ 
for each coloring is calculated by 
$$\left\{
\begin{array}{l}
+f(0,0,1)+f(0,1,2)+f(0,2,0)
=0+2+0=+2, \mbox{ and} \\
+f(0,1,0)+f(0,0,2)+f(0,2,1)
=0+0-2=-2. 
\end{array}\right.$$
Thus we obtain $\Phi_f(D_2,0)=\{\pm2\}$ 
and hence, 
$$\W_f(D_1,\overline{\varphi})-\Phi_f(D_2,0)
=\{-10,-6\}.$$

Next we calculate the set ${\rm Im}(\delta f)$ 
for the map $f$ as above. 
%Recall that 
%\begin{eqnarray*}
%(\delta f)(x,y,z,w) &=& 
%f(x,z,w)-f(x,y,w)+f(x,y,z)\\
%& & -f(x*y,z,w)+f(x*z,y*z,w)-f(x*w,y*w,z*w). 
%\end{eqnarray*}
Since $a*a=a$ holds for any $a\in\Z(3)$, 
we have $(\delta f)(x,y,z,w)=0$ 
for $x=y$, $y=z$, or $z=w$. 
The calculations for 
$x\ne y\ne z\ne w$ are given in Table~$1$. 
For example, 
the calculation for 
$(x,y,z,w)=(0,1,0,1)$ is: 
\begin{eqnarray*}
(\delta f)(0,1,0,1) &=& 
f(0,0,1)-f(0,1,1)+f(0,1,0)\\
& & -f(0*1,0,1)+f(0*0,1*0,1)-f(0*1,1*1,0*1)\\
&=& 
f(0,0,1)-f(0,1,1)+f(0,1,0) \\
& & -f(2,0,1)+f(0,2,1)-f(2,1,2)\\
&=& 0-0+0+2-2+2=+2. 
\end{eqnarray*}

%%%%%%%%%%%%%%
\renewcommand{\arraystretch}{1.1}
\newcommand{\lw}[1]{\smash{
\lower1.8ex\hbox{#1}}}
\begin{table}[htbp]
{\footnotesize 
\begin{tabular}{|c|c|c|c||r|}
\hline
$x$&$y$&$z$&$w$&$\delta f$ \\
\hline
 & &\lw{$0$}&$1$&$+2$\\
\cline{4-5}
 &\lw{$1$}& &$2$&$+7$\\
\cline{3-5}
 & &\lw{$2$}&$0$&$+4$\\
\cline{4-5}
\lw{$0$}& & &$1$&$-1$\\
\cline{2-5}
 & &\lw{$0$}&$1$&$+11$\\
\cline{4-5}
 &\lw{$2$}& &$2$&$+7$\\
\cline{3-5}
 & &\lw{$1$}&$0$&$-4$\\
\cline{4-5}
 & & &$2$&$-8$\\
 \hline
\end{tabular}}
\hspace{10pt}
{\footnotesize 
\begin{tabular}{|c|c|c|c||r|}
\hline
$x$&$y$&$z$&$w$&$\delta f$ \\
\hline
 & &\lw{$1$}&$0$&$+7$\\
\cline{4-5}
 &\lw{$0$}& &$2$&$+5$\\
\cline{3-5}
 & &\lw{$2$}&$0$&$-2$\\
\cline{4-5}
\lw{$1$}& & &$1$&$-4$\\
\cline{2-5}
 & &\lw{$0$}&$1$&$+4$\\
\cline{4-5}
 &\lw{$2$}& &$2$&$-4$\\
\cline{3-5}
 & &\lw{$1$}&$0$&$+1$\\
\cline{4-5}
 & & &$2$&$-7$\\
 \hline
\end{tabular}}
\hspace{10pt}
{\footnotesize 
\begin{tabular}{|c|c|c|c||r|}
\hline
$x$&$y$&$z$&$w$&$\delta f$ \\
\hline
 & &\lw{$1$}&$0$&$+2$\\
\cline{4-5}
 &\lw{$0$}& &$2$&$+4$\\
\cline{3-5}
 & &\lw{$2$}&$0$&$-7$\\
\cline{4-5}
\lw{$2$}& & &$1$&$-5$\\
\cline{2-5}
 & &\lw{$0$}&$1$&$-5$\\
\cline{4-5}
 &\lw{$1$}& &$2$&$-4$\\
\cline{3-5}
 & &\lw{$2$}&$0$&$-1$\\
\cline{4-5}
 & & &$1$&$-2$\\
 \hline
\end{tabular}}
\vspace{10pt}
\label{tab1}
\caption{}
\end{table}
%%%%%%%%%%%%%%

Hence we obtain 
$$\Delta_1(f)=\pm{\rm Im}(\delta f)
=\{0,\pm1,\pm2, \pm4,\pm5,\pm7,\pm 8, \pm11\}.$$
%By the calculations as above, 
{}From the above calculations,
it is easy to check that 
$$\bigl[\W_f(D_1,\overline{\varphi})
-\Phi_f(D_2,0)\bigr]\cap\Delta_i(f)=\emptyset 
\mbox{ for }i=0,1.$$
It follows that $\Omega_3(D_1,D_2)\geq 2$ 
by Proposition~\ref{prop22}. 
The first row of Figure~\ref{fig01} 
shows $\Omega_3(D_1,D_2)\leq 2$, 
which implies $\Omega_3(D_1,D_2)=2$.

\begin{figure}[htb]
\begin{center}
\includegraphics{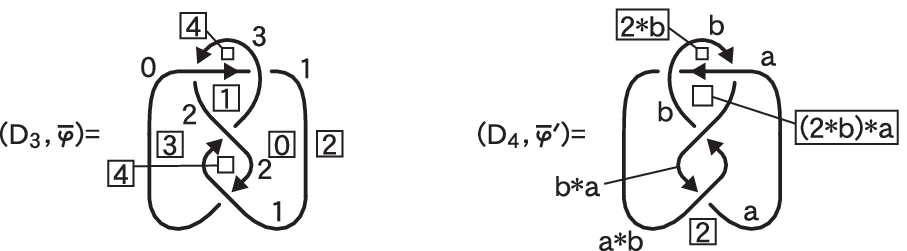}
\caption{}
\label{fig07}
\end{center}
\end{figure}

(ii) We use the function 
$f:\Z(5)^3\rightarrow\Z$ 
defined by 
$f(x,y,z)=(x+y)^3(y+z)(y-z)^3z^5$. %%%
For the extended $5$-coloring 
$\overline{\varphi}=(\varphi,2)$ 
as shown in the left of Figure~\ref{fig07}, 
we have 
\begin{eqnarray*}
\W_f(D_3,\overline{\varphi})
&=& 
-f(3,2,0)-f(0,1,3)+f(4,1,2)+f(4,2,1)\\
&=& 
0+7776-12000+648=-3576.
\end{eqnarray*} 

%%%%%%%%%%%%%%%%%
%%%%%%%%%%%%%%%%%

%Any non-trivial extended $5$-coloring 
%$\overline{\varphi}'=(\varphi',2)$ 
%for $D_4$ is shown in the right of Figure~\ref{fig07}, 
%where $a\ne b\in\Z(5)$. 
All possible non-trivial extended $5$-colorings 
$\overline{\varphi}'=(\varphi',2)$ 
for $D_4$ are indicated in the right of Figure~\ref{fig07},
where $a\ne b\in\Z(5)$. %%%
Then the integer 
$\W_f(D_4,\overline{\varphi}')$ is 
given by 
$$\W_f(D_4,\overline{\varphi}')=
+f(2*b,a,b)+f(2*b,b,a)
-f\bigl((2*b)*a,b,b*a)-f(2,a,a*b).$$
Table~$2$ shows 
all the values $\W_f(D_4,\overline{\varphi}')$ 
for $a\ne b\in\Z(5)$. 
%Hence we obtain 
%%%%%%%%%%%%%%
%$$\renewcommand{\arraystretch}{1.1} 
%\Phi_f(D_4,2)=\left\{
%\begin{array}{rrrrrr}
%-7107048, &  -2814033, & -1919488, & -1299078, & -1269944 \\ 
%-1135296, & -551414, &  -490872, &  -344431, & -173800 \\
 %+142336, &  +207552, &  +326080, &  +587264, &  +889088 \\
 % + 971928, &  +2244931, &   +2937304, &  +3765221, &  +10555072
%\end{array}\right\} . $$
%%%%%%%%%%%%%%

%%%%%%%%%%%%%%
\renewcommand{\arraystretch}{1.1}
\begin{table}[htbp]
\label{tab2}
{\footnotesize 
\begin{tabular}{|c|c||r|}
\hline
$a$&$b$&$\W_f(D_4,\overline{\varphi}')$ \\
\hline
 &$1$&$ +142336$\\
\cline{2-3}
\lw{$0$}&$2$&$ +2244931$\\
\cline{2-3}
 &$3$&$-1269944$\\
\cline{2-3}
 &$4$&$-173800$\\
 \hline
\end{tabular}}
\hspace{5pt}
{\footnotesize 
\begin{tabular}{|c|c||r|}
\hline
$a$&$b$&$\W_f(D_4,\overline{\varphi}')$ \\
\hline
 &$0$&$ +3765221$\\
\cline{2-3}
\lw{$1$}&$2$&$ +207552$\\
\cline{2-3}
 &$3$&$ +587264$\\
\cline{2-3}
 &$4$&$-1299078$\\
 \hline
\end{tabular}}
\hspace{5pt}
{\footnotesize 
\begin{tabular}{|c|c||r|}
\hline
$a$&$b$&$\W_f(D_4,\overline{\varphi}')$ \\
\hline
 &$0$&$ +326080$\\
\cline{2-3}
\lw{$2$}&$1$&$ +971928$\\
\cline{2-3}
 &$3$&$ +2937304$\\
\cline{2-3}
 &$4$&$-1135296$\\
 \hline
\end{tabular}}

\vspace{10pt}
{\footnotesize 
\begin{tabular}{|c|c||r|}
\hline
$a$&$b$&$\W_f(D_4,\overline{\varphi}')$ \\
\hline
 &$0$&$-551414$\\
\cline{2-3}
\lw{$3$}&$1$&$ +889088$\\
\cline{2-3}
 &$2$&$ +10555072$\\
\cline{2-3}
 &$4$&$-344431$\\
 \hline
\end{tabular}}
\hspace{5pt}
{\footnotesize 
\begin{tabular}{|c|c||r|}
\hline
$a$&$b$&$\W_f(D_4,\overline{\varphi}')$ \\
\hline
 &$0$&$-7107048$\\
\cline{2-3}
\lw{$4$}&$1$&$-490872$\\
\cline{2-3}
 &$2$&$-2814033$\\
\cline{2-3}
 &$3$&$-1919488$\\
 \hline
\end{tabular}}
\vspace{10pt}
\caption{}
\end{table}
%%%%%%%%%%%%%%

%On the other hand, 
%we compute the set ${\rm Im}(\delta f)$ 
%by using {\it Maple} and {\it Mathematica} 
%so that it contains $393$ %%%
% integers. 
%We can check that 
On the other hand, we compute the set Im$(\delta f)$ by using
{\it Maple} and {\it Mathematica}, 
and find that it is a set consisting of $393$ integers.
Computer calculations also show that 
$$\bigl[\W_f(D_3,\overline{\varphi})
-\Phi_f(D_4,2)\bigr]\cap\Delta_i(f)=\emptyset 
\mbox{ for } i=0,1,2,$$
and hence we obtain 
$\Omega_3(D_3,D_4)\geq 3$ 
by Proposition~\ref{prop22}. 
The second row of Figure~\ref{fig01} shows 
$\Omega_3(D_3,D_4)\leq 3$, 
which implies $\Omega_3(D_3,D_4)=3$. 

\begin{figure}[htb]
\begin{center}
\includegraphics{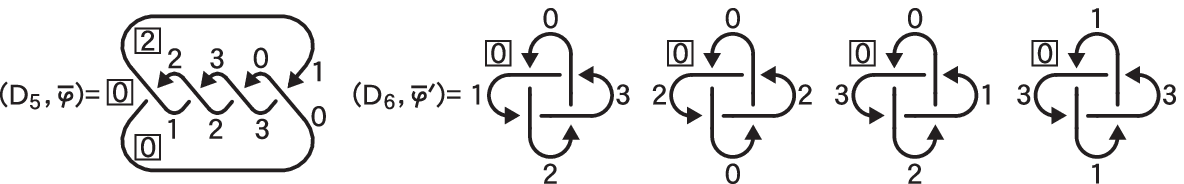}
\caption{}
\label{fig08}
\end{center}
\end{figure}

(iii) We use the function 
$f:\Z(4)^3\rightarrow\Z$ 
defined by 
$f(x,y,z)=(x+y)^2(y-z)^3 z^5$. %%%
Let $\overline{\varphi}=(\varphi,0)$ 
be the extended $4$-coloring for $D_5$ 
as shown in the left of Figure~\ref{fig08}. 
Then it holds that 
\begin{eqnarray*}
\W_f(D_5,\overline{\varphi})
&=& 
+f(2,1,0)+f(2,0,3)+f(2,3,2)+f(2,2,1)\\
&=& 
0-26244+800+16=-25428. %%%
\end{eqnarray*}

%%%%%%%%%%%%%%%%%
%%%%%%%%%%%%%%%%%

Any non-trivial extended $4$-coloring 
$\overline{\varphi}'$ for $D_6$ 
belongs to one of the four cases (up to rotations of $D_6$) 
as shown in the right of Figure~\ref{fig08}. 
The integer $\W_f(D_6,\overline{\varphi}')$ 
for each coloring is calculated by 

$$\left\{
\begin{array}{l}
+f(0,0,1)+f(0,1,2)+f(0,2,3)+f(0,3,0)
=0-32-972+0=-1004, \\
+f(0,0,2)+f(0,2,0)+f(0,0,2)+f(0,2,0)
=2(0+0)=0, \\
+f(0,0,3)+f(0,3,2)+f(0,2,1)+f(0,1,0)
=0+288+4+0=292, \\
+f(0,1,3)+f(0,3,1)+f(0,1,3)+f(0,3,1)
=2(-1944+72)=-3744.
\end{array}\right.$$
Hence we obtain 
$$\Phi_f(D_6,0)=
\{ - 3744, -1004, 0, 292 \}$$
and
$$\W_f(D_5,\overline{\varphi})
-\Phi_f(D_6,0)=
\{ -25720, -25428, -24424, -21684 \}.$$
%%% add:
We note that it is not necessary,  in fact, to check 
the second and the fourth colorings, since the coloring of the left-hand-side
assigns numbers of distinct parities for each component, 
so that this property will be preserved on the right-hand-side. %%%

%On the other hand, 
%we compute the set 
%${\rm Im}(\delta f)$ 
%by using computers as well as (ii) 
%so that it contains $105$ 
% integers. 
%We can check that 
On the other hand, we compute the set Im$(\delta f)$ by using computers again
as we did for the figure-eight knot, and find that it consists of $105$
integers, and also find that
$$\bigl[\W_f(D_5,\overline{\varphi})
-\Phi_f(D_6,0)\bigr]\cap\Delta_i(f)=\emptyset 
\mbox{ for } i=0,1,2,$$
and hence we obtain 
$\Omega_3(D_5,D_6)\geq 3$ 
by Proposition~\ref{prop22}. 
The bottom row of Figure~\ref{fig01} shows 
$\Omega_3(D_5,D_6)\leq 3$, 
which implies $\Omega_3(D_5,D_6)=3$. 
This completes the proof.

\section*{Acknowledgments}

%The first, second, and third authors 
The authors
are 
partially supported by 
NSF Grant DMS $\#0301095$, 
University of South Florida Faculty Development Grant,
NSF Grant DMS $\#0301089$, 
and JSPS Postdoctoral Fellowships for Research Abroad, 
respectively. 
The fourth author expresses his gratitude 
for the hospitality of the University of South Florida 
and the University of South Alabama.

\end{document}